\documentclass[12pt]{article}
\usepackage[latin1]{inputenc}
\usepackage{mathtools}
\usepackage[left=2.5cm, right=2.5cm, top=2cm, bottom=2cm]{geometry}
\usepackage{amsmath,amsthm,amssymb,physics}
\usepackage{url}
\usepackage[colorlinks=true,urlcolor=blue,pdfstartview=FitV, linkcolor=blue,citecolor=blue,bookmarksopen, bookmarksnumbered={true}]{hyperref}

\newtheorem{theorem}{Theorem}[section]
\newtheorem{definition}{Definition}[section]
\newtheorem{example}{Example}[section]

\newtheorem{corollary}{Corollary}[section]
\newtheorem{remark}{Remark}[section]

\begin{document}

\title{\textbf{Infinite Geraghty type extensions and its applications  on integral equations }}

\author{\textbf{Rivu Bardhan$ ^{a}$$^{*}$, Cenap Ozel$ ^{b}$, Liliana Guran$^{c}$}\\ \small{ $^{a}$Tezpur University, Napaam, Tezpur, Assam 784028, India}\\ \small{Email: rivubardhan1@gmail.com}\\ \small{$^b$ Department of Mathematics,	Faculty of Sciences}\\	\small{King Abdulaziz University, Jeddah, Saudi Arabia}.\\ \small{Email: cenap.ozel@gmail.com}\\ \small{$^c$ Department of Pharmaceutical Sciences,}\\	\small{"Vasile Goldi\c s" Western University of Arad, Romania},\\ \small{Email: guran.liliana@uvvg.ro}.\\ \small{$^{*}$Corresponding Author}}

\normalfont
	
\date{}
\maketitle
\begin{abstract}
In this article, we go on to discuss about a series of infinite dimensional extension of the theorems in \cite{Kho}, \cite{F}, \cite{Kg}. We also prove a similar Geraghty type constructions for Fisher (\cite{F}) in infinite dimension, using similar techniques as in  \cite{Kho} and \cite{Kg}. As an application a theory of finding  solutions for the infinite dimensional Fredholm integral equation and Uryshon type integral equation is also provided.
\end{abstract}

\vspace*{5pt}



\section{Introduction and Preliminaries}

\quad The fixed point theory is an essential tool not only in the field of nonlinear analysis but also in almost all applicable branches of importance.  After Banach's contraction principle (\cite{Ban}), over the course of time there were multiple generalizations  have been made for example, Geraghty's Theorem (\cite{G}). Moreover in \cite{Kho} a $k$-dimensional extension of the Geraghty's result is also provided. Never the less, Kannan and Fisher both have given two independent type of contraction operators (see \cite{Kan}, \cite{F}) which are completely independent of the Banach's contraction. Also, in a recent paper by F. Fogh, S. Behnamian and F. Pashaie (see, \cite{Kg}) introduced the concept of Kannan-Geraghty contraction.

In this paper, We extend the $k$-dimensional extension of the Geraghty's Theorem stated in \cite{Kho} to infinite dimension. We also introduce and develop a new notion of $H_k$ contraction for Geraghty-Kannan type operator from \cite{Kg} to provide both $k$ dimensional extension and infinite dimensional extension. We also provide a same family of $k$-dimensional extension for Fisher type contractions and also develop similar infinite dimensional extension in this case. As application of the theorems proved in this paper we define the notion of two independent type of infinite integral equations, namely infinite dimensional Fredholm equation and Uryshon type integral equation and, provide existential proof of their unique solution.

Next, we discuss some of the preliminaries which will be needed later for proving our main theorems.
%Recently  point theory is one of the most useful and essential tools of nonlinear analysis. Banach (see, \cite{Ban}) has given the most important and fundamental theorem of this branch by defining the concept of contraction operators. The following generalization is due to M. Geraghty \cite{G} is one of the many generalizations of it, which has been made over the course of time.

%In this paper, we prove a series of infinite dimensional proper generalizations of results given in \cite{Kho} and \cite{Kg} by finding its fixed points and also develop similar extension for Fisher type mappings.

\begin{theorem}\label{2.1} (\cite{G}) Let $(X, d)$ be a complete metric space and $T : X \to X$ be a
mapping. If $T$ satisfies the following inequality:
\begin{align*}
d(T x, T y) \leq \beta(d(x, y)) d(x, y)
\end{align*}
 where $\beta : [0, \infty) \to [0, 1)$ is a function which satisfies the condition;
 \begin{align*}
\lim\limits_{n\to \infty}\beta(t_n) = 1 ~implies \lim\limits_{n\to \infty}t_n = 0
\end{align*}

Then $T$ has a unique fixed point $u\in X$ and ${T^nx}$ converges to $u$ for each $x \in X.$

\end{theorem}
We denote by $\mathbb{G}$ a set of all functions $\beta$ given in Theorem \ref{2.1}.

\begin{theorem}\label{2.2} (\cite{Kan}) Let $(X, d)$ be a complete metric space and $T : X \to X$ be a
mapping. If $T$ satisfies the following inequality:
\begin{align*}
d(T x, T y) \leq c(d(Tx, x)+ d(Ty, y))
\end{align*}
 where $c\in (0, 1/2)$,

Then $T$ has a unique fixed point $u \in X $.
\end{theorem}

\begin{theorem}\label{2.3} (\cite{F}) Let $(X, d)$ be a complete metric space and $T : X \to X$ be a
mapping. If $T$ satisfies the following inequality:
\begin{align*}
d(T x, T y) \leq c( d(Tx, y)+ d(Ty, x) )
\end{align*}
 where $c\in (0, 1/2)$, $x,y\in X$.

Then $T$ has a unique fixed point $u\in X $.
\end{theorem}

\begin{definition}\label{2.4} (Kannan-Geraghty Map) \cite{Kg} For any complete metric space (X,d) and self map $T:X\to X$, $T$ is defined to be a Kannan-Geraghty self map if and only if it satisfies the following:
		$$d(Tx , Ty)\leq\frac{\beta(d(x,y))}{2}~(d(Tx,x)+d(Ty,y))$$
where, $ \beta \in  \mathbb{G}$, $x,y\in X$.
\end{definition}

\begin{theorem}\label{2.5} (\cite{Kg}) Let $(X, d)$ be a complete metric space and $T : X \to X$ be a
mapping. If $T$ is a Kannan-Geraghty self map on $(X, d)$, then $T$ has a unique fixed point $u \in X $ and for any $x_0 \in X$, $T^nx_0$ converges to $u$.
\end{theorem}

%The above result has so many generalizaions which can be seen in (see [4],[5],[6],[7] etc). \\

\begin{theorem}\label{2.6}(\cite{Kho}) Let $(X, d)$ be a complete metric space and $k \in \mathbb{N} $ let, $T : X^k \to X$
be such that, for all $u_1,u_2, . . . , u_{k+1}\in X$, we have the following inequality,

\begin{align*}
d(T (u_1, . . . , u_k), T (u_2, . . . , u_{k+1}))&\leq
M((u_1, . . . , u_k), (u_2, . . . , u_{k+1})))M((u_1, . . . , u_k), (u_2,. . . ,  u_{k+1}))
\end{align*}
where, $ \beta \in  \mathbb{G}$ and $M : X^k \times X^k \to [0,\infty)$ is as.\\

$M((u_1, . . . , u_k), (u_2,. . . ,  u_{k+1}))
= max\left\lbrace d(u_k, u_{k+1}),~d(u_k, T (u_1, . . . , u_k)),~ d( u_{k+1}, T (u_2,  . . . ,  u_{k+1}))\right\rbrace$.

Then there exists a point $u \in X$ such that $T (u, u, . . . , u) = u.$
\end{theorem}

We denote by
$\mathbb{N}$ (resp. $\mathbb{N}_0$) the set of positive (nonnegative) integers.  The aim of the next section is to generalize and extend Theorem \ref{2.1}, Theorem \ref{2.5} and \ref{2.6} as well as the k-dimensional extension of the result given in \cite{G} to infinite-dimension. We denote the infinite tuples of points $(x_1, x_2, ...)$ by $(x_k)_{i=1} ^{\infty}$ and the infinite tuples $(x_1, x_2, ...x_{k-1}, x_{k},x_{k},...)$ with the $k-$th point repeated by $(x_{j,\hat{k}})_{j=1} ^\infty.$

\section{Main results}
	
\quad In this section we present the main results of the paper. First, let us introduce some new notions concerning the infinite-dimensions of the points. Further let us consider, in all the cases, that $\left(X,d\right)$ is a $complete$ $metric$ $space$ and $k\in\mathbb{N}$.

\begin{definition}\label{d3.1} (Extended $M_k$ function) For any $(u_{i})_{i=1}^{\infty},(v_{i})_{i=1}^{\infty}\in\prod_{i=1}^{\infty} X$ and  $T:\prod_{i=1}^{\infty} X\to X$, we define the extended function $M_{k}:\prod_{l=1}^\infty X \times \prod_{l=1}^\infty X\rightarrow X$ as follows,

$$M_{k}\left((u_{i})_{i=1}^{\infty},(v_{i})_{i=1}^{\infty}\right)=\left\{\begin{array}{ll} {\max\left\lbrace\sup\limits_{l\geqslant k} d\left(u_{l},v_{l}\right),\sup\limits_{l\geqslant k} d\left(u_{l},T(u_{i,\hat{l}})_{i=1}^{\infty}\right), \sup\limits_{l\geqslant k} d\left(v_{l},T(v_{i,\hat{l}})_{i=1}^{\infty}\right)\right\rbrace},\\ \hspace{20pt}{\text{ if all supremum exist}}\\
{\max\left\lbrace d\left(u_{k},v_{k}\right), d\left(u_{k},T(u_{i,\hat{k}})_{i=1}^{\infty}\right), d\left(v_{k},T(v_{i,\hat{k}})_{i=1}^{\infty}\right)\right\rbrace},\\ \hspace{20pt} {\text{ if one fails to exists}}\\
\end{array}\right.$$

\end{definition}

\begin{example}
 For $k=1\text{, }X=[0,1]$ the operator $T:\prod_{i=1}{\infty}X\to X$ defined as $T((u_{i})_{i=1}^{\infty})=u_k$, $M_{k}\left((u_{i})_{i=1}^{\infty},(v_{i})_{i=1}^{\infty}\right)= \abs{u_k-v_k}=\abs{u_1-v_1}$
\end{example}

\begin{definition}\label{d3.2} ($H_k$ contraction) An operator $T:\prod_{l=1}^\infty X\rightarrow X$ is called a $H_{k}$ $contraction$ if and only if it satisfies the following inequality,
\begin{equation}\label{hc}
		d\left(T\left((u_{i,\widehat{k}})_{i=1}^{\infty}\right), T\left((v_{i,\widehat{k}})_{i=1}^{\infty}\right)\right) \leqslant\beta\left(M_{k}\left((u_{i,\widehat{k}})_{i=1}^{\infty},(v_{i,\widehat{k}})_{i=1}^{\infty}\right)\right)M_{k}\left((u_{i,\widehat{k}})_{i=1}^{\infty},(v_{i,\widehat{k}})_{i=1}^{\infty}\right),
\end{equation}
for all $u_1,~u_2,~\ldots,~u_{k},~ v_1,~ v_2, \ldots,~ v_k\in X$, where $\beta\in\mathbb{G}$.
\end{definition}

\begin{example}
 For example the operator $T:\prod_{i=1}^{\infty}X\to X$ where $X=[0,1]$, defined by $T((x_i)_{i=1}^{\infty}):=cx_k$ is a $H_k$ contraction where $k\in\mathbb{N}$ be any fixed number and $c\in(0,1)$.
\end{example}

\begin{definition}\label{d3.3} (Kannan-Geraghty $H_k$ contraction) An operator $T:\prod_{l=1}^{k} X\rightarrow X$ is called a Kannan-Geraghty contraction of dimension $k$ if and only if it satisfies the following inequality
\begin{equation*}
  d\left(T\left((u_{i})_{i=1}^{k}\right), T\left((v_{i})_{i=1}^{k}\right)\right) \leqslant\frac{\beta(d(u_k,v_{k}))}{2}(d(T\left((u_{i})_{i=1}^{k}\right), u_k)+ d(T\left((v_{i})_{i=1}^{k}\right), v_{k})),
\end{equation*}		
for all $u_1, u_2,\ldots,u_{k}, v_1, v_2, \ldots, v_k\in X$, where $\beta\in\mathbb{G}$.
\end{definition}

\begin{definition}\label{d3.4} (Extended Kannan-Geraghty $H_k$ contraction) An operator $T:\prod_{l=1}^{\infty} X\rightarrow X$ is called a Extended Kannan-Geraghty $H_k$ contraction if and only if it satisfies the following inequality
\begin{equation*}
		d\left(T\left((u_{i,\widehat{k}})_{i=1}^{\infty}\right),T\left((v_{i,\widehat{k}})_{i=1}^{\infty}\right)\right) \leqslant\frac{\beta\left(d(u_k,v_{k})\right)}{2}(d(T\left((u_{i,\widehat{k}})_{i=1}^{\infty}\right),u_k)+d(T\left((v_{i,\widehat{k}})_{i=1}^{\infty}\right),v_{k})),
\end{equation*}		
for all $u_1,~ u_2,\ldots,~u_{k}, v_1,~ v_2, \ldots,~ v_k\in X$, where $\beta\in\mathbb{G}$.

\end{definition}

\begin{definition}\label{d3.5} (Fisher-Geraghty $H_k$ contraction)  An operator $T:\prod_{l=1}^{k} X\rightarrow X$ is called a Fisher-Geraghty $H_k$ contraction of dimension $k$ if and only if it satisfies the following inequality,
\begin{equation*}
	d\left(T\left((u_{i})_{i=1}^{k}\right), T\left((v_{i})_{i=1}^{k}\right)\right) \leqslant\frac{\beta(d(u_k,v_k))}{2}(d(T\left((u_{i})_{i=1}^{k}\right),v_{k})+d(T\left((v_{i})_{i=1}^{k}\right),u_{k})),
\end{equation*}
for all $u_1,~ u_2,~\ldots,~u_{k},~ v_1,~ v_2,~ \ldots,~ v_k\in X$, where $\beta\in\mathbb{G}$.
\end{definition}

\begin{definition}\label{d3.6} (Extended Fisher-Geraghty $H_k$ contraction)  An operator $T:\prod_{l=1}^{\infty} X\rightarrow X$ is called a Extended Kannan-Geraghty $H_k$ contraction if and only if it satisfies the following inequality
\begin{equation*}
d\left(T\left((u_{i,\widehat{k}})_{i=1}^{\infty}\right),T\left((v_{i,\widehat{k}})_{i=1}^{\infty}\right)\right) \leqslant\frac{\beta\left(d(u_k,v_{k})\right)}{2}(d(T\left((u_{i,\widehat{k}})_{i=1}^{\infty}\right),v_{k})+d(T\left((v_{i,\widehat{k}})_{i=1}^{\infty}\right),u_{k})),
\end{equation*}
for all $u_1,~ u_2,~\ldots,~u_{k},~ v_1,~ v_2,~ \ldots,~ v_k\in X$, where $\beta\in\mathbb{G}$.
\end{definition}

\begin{example}
 As an example to Definition \ref{d3.3} and Definition \ref{d3.5} we can consider for $X=[0,1]$, $T:\prod_{i=1}^{k}X\to~X$ defined as $T((x_i)_{i=1}^{k})=cx_k$ for some $c\in(0,\frac{1}{2})$.
\end{example}

\begin{example}
 As an example to Definition \ref{d3.4} and Definition \ref{d3.6} we can consider for $X=[0,1]$, $T:\prod_{i=1}^{\infty} X\to~X$ defined as $T((x_i)_{i=1}^{\infty})=cx_k$ for some $c\in(0,\frac{1}{2})$.
\end{example}
\noindent Let $T:\prod_{i=1}^{k}X\to~X$ be any operator. Let us choose $x_1, x_2, \ldots, x_k\in X$ Now we define the following:
\begin{definition}\label{d3.7}($k$-Picard sequence with respect to the operator $T$) %\textcolor[rgb]{1.00,0.00,0.00}{Please add the Reference item where we find this operator and cite it here, after definition}
The $k$-Picard sequence with respect to the operator $T$ based on the base point set $\lbrace x_1, x_2, \ldots, x_k\rbrace$ is definded as  $x_{n+k}:=T((x_{n+i-1})_{i=1}^{k})$, for all $n\geq k$.
\end{definition}

\begin{example}
 If we fix $k=1$ then the base point set is singleton and the $1$- Picard sequence with respect to $T$ based on $\lbrace x_0\rbrace$ is basically the Picard sequence of $T$ based on the base point $\lbrace x_0\rbrace$ defined by $x_n:= T(x_{n-1}) ~\forall n\geqslant 1$ for some $x_0\in X$.
\end{example}

\noindent Let $T:\prod_{i=1}^{\infty}X\to~X$ be any operator. Let us choose $x_1, x_2, \ldots, x_k\in X$ Now we define the following:

\begin{definition}\label{d3.8}(Infinite $k$-Picard sequence with respect to the operator $T$) %\textcolor[rgb]{1.00,0.00,0.00}{Please add the reference item where we find this operator and cite it here, after definition}
The $k$-Picard sequence with respect to the operator $T$ based on the base point set $\lbrace x_1, x_2, \ldots, x_k\rbrace$ is definded as  $x_{n+k}:=T((x_{n+i-1,\widehat{n+k-1}})_{i=1}^{\infty})$, for all $n\geq k$.
\end{definition}

\begin{example}
If we fix $k=1$ then the base point set is singleton and the infinite $1$- Picard sequence with respect to $T$ based on the base point $\lbrace x_0\rbrace$ is basically the sequence  defined by $x_n:= T((x_{n-1})_{n=1}^{\infty}) ~\forall n\geqslant 1$ for some $x_0\in X$.
\end{example}

Let us give our first main fixed point result which is a generalization of Banach contraction principle with respect to the infinite-dimensional notion introduced in our paper.
\begin{theorem}\label{3.1}
Let $\left(X,d\right)$ be a complete metric space. $T:\prod_{l=1}^\infty X \rightarrow X$ is a $ H_{k}$ contraction for some $k\in\mathbb{N}$ . Then there exists $u\in X$ such that $T((u)_{i=1}^{\infty})=u$ and the infinite $k$-Picard sequence for $T$ converges to $u$.
\end{theorem}
\begin{proof}
		Let $x_1,x_2\ldots,x_k\in X$. Then infinite $k$-Picard sequence is defined as follows:
		\newline
		$\forall~n\in\mathbb{N} $ we define,
		\begin{equation}\label{s}
		x_{n+k}:=T\left((x_{n+i-1,\widehat{n+k-1}})_{i=1}^{\infty}\right)
		\end{equation}
		\newline	
	
			We claim that $\left\lbrace d\left(x_{n+k},x_{n+k+1}\right)\right\rbrace_{n\in\mathbb{N}}$ is convergent.

		From Definition (\ref{d3.2})
		\begin{flushleft}
			$d\left(x_{n+k+1},x_{n+k+2}\right)=d\left(T\left((x_{n+i,\widehat{n+k}})_{i=1}^{\infty}\right), T\left((x_{n+i+1,\widehat{n+1+k}})_{i=1}^{\infty}\right)\right)$
			
			$\leqslant\beta\left(M_{k}\left((x_{n+i,\widehat{n+k}})_{i=1}^{\infty}, (x_{n+1+i,\widehat{n+1+k}})_{i=1}^{\infty}\right)\right)M_{k}\left((x_{n+i,\widehat{n+k}})_{i=1}^{\infty}, (x_{n+1+i,\widehat{n+1+k}})_{i=1}^{\infty}\right).$
		
		\end{flushleft}
		\begin{flushleft}
			Let
$M_{k}\left((x_{n+i,\widehat{n+k}})_{i=1}^{\infty},(x_{n+i+1,\widehat{n+1+k}})_{i=1}^{\infty}\right)=\max\left\lbrace\sup\limits_{l\geqslant k} d\left(u_{l},v_{l}\right),\sup\limits_{l\geqslant k} d\left(u_{l},T(u_{i,\widehat{l}})_{i=1}^{\infty}\right),\sup\limits_{l\geqslant k} d\left(v_{l},T(v_{i,\widehat{l}})_{i=1}^{\infty}\right)\right\rbrace$
		\end{flushleft}
		where
		\begin{gather}
		u_{l}:=x_{n+l},\: \forall\: 1\leqslant l\leqslant k
		\\
		u_{l}:=x_{n+k},\: \forall\: l> k
		\\
		v_{l}:=x_{n+l+1},\: \forall\: 1\leqslant l\leqslant k
		\\
		u_{l}:=x_{n+k+1},\: \forall\: l> k.
		\end{gather}
		But
		\begin{flushleft}
			$\sup\limits_{l\geqslant k} d\left(u_{l},v_{l}\right)=d\left(x_{n+k},x_{n+k+1}\right)$,\\
			$\sup\limits_{l\geqslant k} d\left(u_{l},T\left((u_{i,\widehat{l}})_{i=1}^{\infty}\right)\right)=d\left(x_{n+k},T\left((x_{n+i,\widehat{n+k}})_{i=1}^{\infty}\right)\right)$,\\
			$\sup\limits_{l\geqslant k} d\left(v_{l},T\left((v_{i,\widehat{l}})_{i=1}^{\infty}\right)\right)=d\left(x_{n+k+1},T\left((x_{n+i+1,\widehat{n+k+1}})_{i=1}^{\infty}\right)\right)$.
		\end{flushleft}
		\begin{flushleft}
			$\implies M_{k}\left((x_{n+i,\widehat{n+k}})_{i=1}^{\infty},(x_{n+i+1,\widehat{n+k+1}})_{i=1}^{\infty}\right)$\\$=\max\left\lbrace d\left(x_{n+k},x_{n+k+1}\right), d\left(x_{n+k},T\left((x_{n+i,\widehat{n+k}})_{i=1}^{\infty}\right)\right),d\left(x_{n+k+1},T\left((x_{n+1+i,\widehat{n+k+1}})_{i=1}^{\infty}\right)\right)\right\rbrace$
			$\implies M_{k}\left((x_{n+i,\widehat{n+k}})_{i=1}^{\infty},(x_{n+1+i,\widehat{n+k+1}})_{i=1}^{\infty}\right)=\max\left\lbrace d\left(x_{n+k},x_{n+k+1}\right), d\left(x_{n+k},x_{n+k+1}\right),d\left(x_{n+k+1},x_{n+k+2}\right)\right\rbrace$
			$=\max\left\lbrace d\left(x_{n+k},x_{n+k+1}\right),d\left(x_{n+k+1},x_{n+k+2}\right)\right\rbrace$.
		\end{flushleft}

		If $M_{k}\left((x_{n+i,\widehat{n+k}})_{i=1}^{\infty},(x_{n+i+1,\widehat{n+k+1}})_{i=1}^{\infty}\right)=d\left(x_{n+k+1},x_{n+k+2}\right)$ then, from Definition (\ref{d3.2}), we get
		
				\begin{flushleft}
			$d\left(x_{n+k+1},x_{n+k+2}\right)=d\left(T\left((x_{n+i,\widehat{n+k}})_{i=1}^{\infty}\right),T\left((x_{n+i+1,\widehat{n+k+1}})_{i=1}^{\infty}\right)\right)$
			\\
			$\leqslant\beta\left(M_{k}\left((x_{n+i,\widehat{n+k}})_{i=1}^{\infty},(x_{n+i+1,\widehat{n+k+1}})_{i=1}^{\infty}\right)\right)M_{k}\left((x_{n+i,\widehat{n+k}})_{i=1}^{\infty},(x_{n+i+1,\widehat{n+k+1}})_{i=1}^{\infty}\right)$
			\\
			$=\beta\left(d\left(x_{n+k+1},x_{n+k+2}\right)\right)d\left(x_{n+k+1},x_{n+k+2}\right)$
			\\
			$<d\left(x_{n+k+1},x_{n+k+2}\right).$
		\end{flushleft}

Then we get a contradiction.
	
			Hence $M_{k}\left((x_{n+i,\widehat{n+k}})_{i=1}^{\infty},(x_{n+i+1,\widehat{n+k+1}})_{i=1}^{\infty}\right)=d\left(x_{n+k},x_{n+k+1}\right).$
		\begin{flushleft}
			$\implies d\left(x_{n+k+1},x_{n+k+2}\right)\leqslant\beta\left(d\left(x_{n+k},x_{n+k+1}\right)\right)d\left(x_{n+k},x_{n+k+1}\right)<d\left(x_{n+k},x_{n+k+1}\right)$
		\end{flushleft}
		$\implies \left\lbrace x_{n+k}\right\rbrace_{n\in\mathbb{N}}$ is a strictly decreasing sequence.\\
		
	Using Monotone Convergence Theorem, for some $(r\geqslant 0)\in \mathbb{R}$, we get
		\begin{equation}\label{l}
		\lim_{n\to\infty}d( x_{n+k},x_{n+k+1})=r.\\
		\end{equation}

		We claim $\lim\limits_{n\to\infty}d( x_{n+k},x_{n+k+1})=0$.
	
		Suppose the contrary and assuming that $r>0$ we get
		\begin{flushleft}
			$$d\left(x_{n+k+1},x_{n+k+2}\right)\leqslant\beta\left(d\left(x_{n+k},x_{n+k+1}\right)\right)d\left(x_{n+k},x_{n+k+1}\right)$$ $$\implies\frac{d\left(x_{n+k+1},x_{n+k+2}\right)}{d\left(x_{n+k},x_{n+k+1}\right)}\leqslant\beta\left(d\left(x_{n+k},x_{n+k+1}\right)\right).$$
		\par This implies $\lim\limits_{n\to\infty}\beta\left(d\left(x_{n+k},x_{n+k+1}\right)\right)\geqslant 1.$
		\par Since
			$\beta\in\mathbb{G}\text{ implies}\lim\limits_{n\to\infty}\beta\left(d\left(x_{n+k},x_{n+k+1}\right)\right)\leqslant 1$
			$\implies\lim\limits_{n\to\infty}\beta\left(d\left(x_{n+k},x_{n+k+1}\right)\right)=1$
we obtain 	
			\begin{equation}
			\lim\limits_{n\to\infty}d\left(x_{n+k},x_{n+k+1}\right)=0.\label{8}
			\end{equation}
		\end{flushleft}
		
		We now claim that $\left\lbrace x_{n+k}\right\rbrace_{n\in\mathbb{N}}$ is Cauchy and we prove it by contradiction. If we suppose the contrary $\exists~\varepsilon>0$ such that we can find some subsequences $\left\lbrace x_{\left(m\right)+k}\right\rbrace_{p\in\mathbb{N}}$, $\left\lbrace x_{n\left(p\right)+k}\right\rbrace_{p\in\mathbb{N}}$ with $m\left(p\right)>n\left(p\right)>p$ such that for every $p$ we have
			\begin{equation}
			d\left(x_{m\left(p\right)+k},x_{n\left(p\right)+k}\right)\geqslant\varepsilon\label{9}
			\end{equation}
			
			Moreover, corresponding to each $n(p)$ we can choose least of such $m\left(p\right)$ satisfying (\ref{9}). Then
			\begin{equation}\label{10}
			d\left(x_{m\left(p\right)+k-1},x_{n\left(p\right)+k}\right)<\varepsilon.
			\end{equation}
			
			From (\ref{8}),(\ref{10}) and using triangle inequality we get
			\begin{equation}
			d\left(x_{m\left(p\right)+k-1},x_{n\left(p\right)+k-1}\right)\leqslant d\left(x_{m\left(p\right)+k-1},x_{n\left(p\right)+k}\right)+d\left(x_{n\left(p\right)+k-1},x_{n\left(p\right)+k}\right)\label{11}
			\end{equation}
			\begin{center}
				$<\varepsilon+d\left(x_{n\left(p\right)+k-1},x_{n\left(p\right)+k}\right)$
			\end{center}
			and			
			\begin{equation}\label{12}
			\varepsilon \leqslant d\left(x_{n\left(p\right)+k},x_{m\left(p\right)+k}\right)\\
			~\leqslant d\left(x_{n\left(p\right)+k},x_{n\left(p\right)+k-1}\right)+d\left(x_{n\left(p\right)+k-1},x_{m\left(p\right)+k-1}\right)+d\left(x_{m\left(p\right)+k-1},x_{m\left(p\right)+k}\right).
			\end{equation}
			
			If $p\to\infty$ in (\ref{11}) and using (\ref{12}) we get
			\begin{equation}\label{13}
			\lim_{p\to\infty}d\left(x_{m\left(p\right)+k-1},x_{n\left(p\right)+k-1}\right)=\varepsilon.
			\end{equation}
	
		\begin{flushleft}
			On the other hand, if
			$M_{k}\left((x_{n(p)+i-1,\widehat{n(p)+k-1}})_{i=1}^{\infty},(x_{m(p)+i-1,\widehat{m(p)+k-1}})_{i=1}^{\infty}\right)$ \\$=\max\left\lbrace\sup\limits_{l\geqslant k-1} d\left(u_{l},v_{l}\right),\sup\limits_{l\geqslant k-1} d\left(u_{l},T\left((u_{i,\widehat{l}})_{i=1}^{\infty}\right)\right),\sup\limits_{l\geqslant k-1} d\left(v_{l},T\left((v_{i,\widehat{l}})_{i=1}^{\infty}\right)\right)\right\rbrace$
			where
			\begin{gather}
			u_{l}:=x_{n\left(p\right)+l-1},\: \forall\: 1\leqslant l\leqslant k
			\\
			u_{l}:=x_{n\left(p\right)+k-1},\: \forall\: l> k
			\\
			v_{l}:=x_{m\left(p\right)+l-1},\: \forall\: 1\leqslant l\leqslant k
			\\
			u_{l}:=x_{m\left(p\right)+k-1},\: \forall\: l> k
			\end{gather}
			$\implies M_{k}\left((x_{n(p)+i-1,\widehat{n(p)+k-1}})_{i=1}^{\infty},(x_{m(p)+i-1,\widehat{m(p)+k-1}})_{i=1}^{\infty}\right)$
			
			$=\max\lbrace d\left(x_{n(p)+k-1},x_{m(p)+k-1}\right),  d\left(x_{n(p)+k-1},T\left((x_{n(p)+i-1,\widehat{n(p)+k-1}})_{i=1}^{\infty}\right)\right)$,\\
			$\hspace{6cm} d\left(x_{m\left(p\right)+k-1},T\left((x_{m(p)+i-1,\widehat{m(p)+k-1}})_{i=1}^{\infty})\right)\right)\rbrace$\\
			$=\max\lbrace d\left(x_{n(p)+k-1},x_{m(p)+k-1}\right), d(x_{n(p)+k-1},x_{n(p)+k}), d(x_{m(p)+k-1},x_{m(p)+k})\rbrace $.
		\end{flushleft}
		Using (\ref{8}) and (\ref{13}) we get
		\begin{equation}
		\lim_{p\to\infty}M_{k}\left((x_{n(p)+i-1,\widehat{n(p)+k-1}})_{i=1}^{\infty},(x_{m(p)+i-1,\widehat{m(p)+k-1}})_{i=1}^{\infty}\right)=\varepsilon.\label{18}
		\end{equation}
		By (\ref{hc}) and (\ref{18}) we get
		\begin{align}\label{e8}
&\varepsilon\leq d\left(x_{n\left(p\right)+k},x_{m\left(p\right)+k}\right)=d\left(T(x_{n(p)+i-1,\widehat{n(p)+k-1}})_{i=1}^{\infty},T(x_{m(p)+i-1,\widehat{m(p)+k-1}})_{i=1}^{\infty})\right)\nonumber\\ &\leq\beta\left(M_{k}\left((x_{n(p)+i-1,\widehat{n(p)+k-1}})_{i=1}^{\infty},(x_{m(p)+i-1,\widehat{m(p)+k-1}})_{i=1}^{\infty}\right)\right)\\
&\ \ \ \ \ M_{k}\left((x_{n(p)+i-1,\widehat{n(p)+k-1}})_{i=0}^{\infty},(x_{m(p)+i-1,\widehat{m(p)+k-1}})_{i=1}^{\infty}\right)\nonumber\\
&<M_{k}\left((x_{n(p)+i-1,\widehat{n(p)+k-1}})_{i=1}^{\infty},(x_{m(p)+i-1,\widehat{m(p)+k-1}})_{i=1}^{\infty}\right).\nonumber
			\end{align}

			If we suppose
			    $$\sup\limits_{l\geqslant k} d\left(u_{l},v_{l}\right)=d\left(x_{n(p)+k-1},x_{m(p)+k-1}\right)$$
				$$\sup\limits_{l\geqslant k} d\left(u_{l},T(u_{i,\widehat{l}})_{i=1}^{\infty}\right)=d\left(x_{n(p)+k-1},T\left((x_{n(p)+i-1,\widehat{n(p)+k-1}})_{i=1}^{\infty}\right)\right)$$
				$$\sup\limits_{l\geqslant k} d\left(v_{l},T(v_{i,\widehat{l}})_{i=1}^{\infty}\right)=d\left(x_{m(p)+k-1},T\left((x_{m(p)+i-1,\widehat{m(p)+k-1}})_{i=1}^{\infty}\right)\right)$$
\begin{align*}
			&\implies M_{k}\left((x_{n(p)+i-1,\widehat{n(p)+k-1}})_{i=1}^{\infty},(x_{m(p)+i-1,\widehat{m(p)+k-1}})_{i=1}^{\infty}\right)\\&=\max\lbrace d\left(x_{n(p)+k},x_{m(p)+k}\right),d\left(x_{n(p)+k-1}, T\left((x_{n(p)+i-1,\widehat{n(p)+k-1}})_{i=1}^{\infty}\right)\right),\\ &d\left(x_{m\left(p\right)+k-1},T\left((x_{m(p)+i-1,\widehat{m(p)+k-1}})_{i=1}^{\infty})\right)\right)\rbrace
			\\&=\max\lbrace d\left(x_{n(p)+k}, x_{m(p)+k}\right), d(x_{n(p)+k-1}, x_{n(p)+k}), d(x_{m(p)+k-1}, x_{m(p)+k})\rbrace
	\end{align*}	
		\begin{flushleft}
			Letting $p\to\infty$ and using (\ref{hc}), (\ref{e8}), (\ref{8}) and (\ref{18}) we get\\ $\lim\limits_{p\to\infty}\beta\left(M_{k}\left((x_{n(p)+i-1,\widehat{n(p)+k-1}})_{i=1}^{\infty},(x_{m(p)+i-1,\widehat{m(p)+k-1}})_{i=1}^{\infty}\right)\right)$\\$\lim\limits_{p\to\infty}M_{k}\left((x_{n(p)+i-1,\widehat{n(p)+k-1}})_{i=1}^{\infty},(x_{m(p)+i-1,\widehat{m(p)+k-1}})_{i=1}^{\infty}\right)\geqslant\varepsilon$.\\
	
			From (\ref{18})		$$\beta\left(M_{k}\left((x_{n(p)+i-1,\widehat{n(p)+k-1}})_{i=1}^{\infty},(x_{m(p)+i-1,\widehat{m(p)+k-1}})_{i=1}^{\infty}\right)\right)\geqslant 1$$
			Since $\beta\in\mathbb{G}$.
			$\implies \lim\limits_{p\to\infty}\beta\left(M_{k}\left((x_{n(p)+i-1,\widehat{n(p)+k-1}})_{i=1}^{\infty},(x_{m(p)+i-1,\widehat{m(p)+k-1}})_{i=1}^{\infty}\right)\right)= 1$
		$\implies\lim\limits_{p\to\infty}M_{k}\left((x_{n(p)+i-1,\widehat{n(p)+k-1}})_{i=1}^{\infty},(x_{m(p)+i-1,\widehat{m(p)+k-1}})_{i=1}^{\infty}\right)=0$
		\end{flushleft}
		We get a contradiction to (\ref{18}).
		Hence $\{x_{n+k}\}_{n\in\mathbb{N}}$ is a Cauchy sequence.
		
		Since $\left(X,d\right)$ is complete there exists $u\in X$ such that
		\begin{equation}\label{e9}
		\lim\limits_{n\to\infty}x_{n+k}=u.
		\end{equation}

		Then we claim that $d\left(T\left((u)_{i=1}^{\infty}\right),u\right)=0$.
	
		\begin{flushleft}
			If we suppose contrary, then $d\left(u,T\left((u)_{i=1}^{\infty}\right)\right)>0$,
			\begin{align*}
			M\left((x_{n+i-1,\widehat{n+k-1}})_{i=1}^{\infty},\left((u)_{i=1}^{\infty}\right)\right)&=\max\left\lbrace d(x_{n+k-1},u),~d\left(x_{n+k-1},T(x_{n+i-1,\widehat{n+k-1}})_{i=1}^{\infty}\right),d\left(u,T\left((u)_{i=1}^{\infty}\right)\right)\right\rbrace\\
		&=\max\left\lbrace d\left(x_{n+k-1},u\right),d\left(x_{n+k},x_{n+k-1}\right),d\left(u,T((u)_{i=1}^{\infty})\right)\right\rbrace.\\
			\end{align*}
		
		\end{flushleft}
		Letting $n\to\infty$ and using (\ref{8}), (\ref{e9}) we get	$$\lim\limits_{n\to\infty}M\left((x_{n+i-1,\widehat{n+k-1}})_{i=1}^{\infty},(u)_{i=1}^{\infty}\right)=d\left(u,T((u)_{i=1}^{\infty})\right)\neq 0.$$
		\begin{flushleft}
			Then we have
			\begin{align*}
d\left(u,T\left((u)_{i=1}^{\infty}\right)\right)&\leqslant d\left(x_{n+k},u\right)+d\left(x_{n+k},T((u)_{i=1}^{\infty})\right)\\ &=d\left(x_{n+k},u\right)+d\left(T\left((x_{n+i-1,\widehat{n+k-1}})_{i=1}^{\infty}\right),T((u)_{i=1}^{\infty})\right)
			\\ & \leqslant d\left(x_{n+k},u\right)+\beta\left(M\left(
			(x_{n+i-1,\widehat{n+k-1}})_{i=1}^{\infty},(u)_{i=1}^{\infty}\right)\right)M\left(
			(x_{n+i-1,\widehat{n+k-1}})_{i=1}^{\infty},(u)_{i=1}^{\infty}\right)
			\\
			& \leqslant d\left(x_{n+k},u\right)+\beta\left(d\left(u,T((u)_{i=1}^{\infty})\right)\right)d\left(u,T((u)_{i=1}^{\infty})\right).
			\end{align*}
			\end{flushleft}

		Letting $n\to\infty$ and using (\ref{e9}) we get $\lim_{n\to\infty}\beta\left(d\left(u,T((u)_{i=1}^{\infty})\right)\right)\geqslant 1$.

Since $\beta\in\mathbb{G}\implies\lim\limits_{n\to\infty}\beta\left(d\left(u,T((u)_{i=1}^{\infty})\right)\right)=1.$

Then we have $\lim\limits_{n\to\infty} d\left(u,T((u)_{i=1}^{\infty})\right)=0$, which implies $d\left(u,T((u)_{i=1}^{\infty})\right)=0$. So we get a contradiction.
		
Hence, $d\left(u,T((u)_{i=1}^{\infty})\right)=0\implies T((u)_{i=1}^{\infty}))=u$.
\end{proof}
		
\begin{remark} Theorem \ref{3.1} is a proper generalization to the Theorem \ref{2.6} since in case of the simplest operator on $\prod_{i=1}^{\infty}X\to X$ the contraction condition of Theorem \ref{2.6} is not applicable but on the other hand the $H_k$ contraction (see Definition \ref{d3.2} ) is easily applicable for the infinite case. Also if we restrict the operator to any finite $k$ dimension through an easy calculation it is obvious that it is an equivalent statement of Theorem \ref{2.6}.

\end{remark}

\begin{theorem}\label{3.2}.
Let $\left(X,d\right)$ be a complete metric space and $T:\prod_{l=1}^{\infty} X \rightarrow X$ be an extended Kannan-Geraghty $H_k$ contraction for some $k\in\mathbb{N}$. Then $\exists u\in X$ such that $T((u)_{i=1}^{\infty})=u$ and for any $x_1, \ldots, x_k \in X$, the infinite $k$-Picard sequence converges to $u$.\\
\end{theorem}
\begin{proof}
  Let $x_1, \ldots, x_k\in X$. For all $n\in\mathbb{N}$ we define the infinite $k$-Picard sequence as follows
		\begin{equation}\label{s}
		x_{n+k}:=T\left((x_{n+i-1,\widehat{n+k-1}})_{i=1}^{\infty}\right).
		\end{equation}

		We claim that $\lim_{n\to\infty}d(x_{n+k},x_{n+k+1})=0$. Then we have
		\begin{flushleft}			$d\left(x_{n+k+1},x_{n+k+2}\right)=d\left(T\left((x_{n+i,\widehat{n+k}})_{i=1}^{\infty}\right),T\left((x_{n+i+1,\widehat{n+1+k}})_{i=1}^{\infty}\right)\right)$
			\\		$\leqslant\frac{\beta\left(d(x_{n+k},x_{n+k+1})\right)}{2}(d(T\left((x_{n+i,\widehat{n+k}})_{i=1}^{\infty}\right),x_{n+k})+d(T\left((x_{n+i+1,\widehat{n+k+1}})_{i=1}^{\infty}\right),x_{n+k+1}))$
			\\
			$<~\frac{1}{2}(d(x_{n+k},x_{n+k+1})+d(x_{n+k+1},x_{n+k+2}))\implies d(x_{n+k+1},x_{n+k+2})<d(x_{n+k},x_{n+k+1})$.
		\end{flushleft}

Then our sequence  is a monotone decreasing sequence which is bounded below. Then there exists $r\geq~0$ such that
\begin{equation}\label{1}
  \lim\limits_{n\to~0}d(x_{n+k},x_{n+k+1})=r.
\end{equation}

We claim that $r=0$.

If we consider the contrary and suppose $r>0$ we have
\begin{flushleft}
  $\lim\limits_{n\to\infty}d(x_{n+k+1},x_{n+k+2})\leqslant\frac{\lim\limits_{n\to\infty}\beta\left(d(x_{n+k},x_{n+k+1})\right)}{2}(\lim\limits_{n\to\infty}d(x_{n+k},x_{n+k+1})+\lim\limits_{n\to\infty}d(x_{n+k+1},x_{n+k+2}))$

  $\implies\frac{2\lim\limits_{n\to\infty}d(x_{n+k+1},x_{n+k+2})}{\lim\limits_{n\to\infty}d(x_{n+k},x_{n+k+1})+\lim\limits_{n\to\infty}d(x_{n+k+1},x_{n+k+2})}\leqslant~\lim\limits_{n\to\infty}\beta\left(d(x_{n+k},x_{n+k+1})\right)$.

  \par From (\ref{1}) we get

  $\frac{2r}{2r}\leqslant~\beta\left(d(x_{n+k},x_{n+k+1})\right)\implies\lim\limits_{n\to\infty}\beta\left(d(x_{n+k},x_{n+k+1})\right)\geqslant~1$.
\end{flushleft}
		
Since $\beta\in\mathbb{G}$ then $\lim\limits_{n\to\infty}\beta\left(d(x_{n+k},x_{n+k+1})\right)\leqslant~1$.

Using the well known "Sandwich Theorem" we obtain
\begin{equation}
  \lim\limits_{n\to\infty}\beta\left(d(x_{n+k},x_{n+k+1})\right)=1\implies\lim\limits_{n\to\infty}d(x_{n+k},x_{n+k+1})=0.
\end{equation}

\begin{flushleft}

 Further, we have

   $d(x_{n+k},x_{m+k})=d(T((x_{n+i-1,\widehat{n+k-1}})_{i=1}^{\infty}),T((x_{m+i-1,\widehat{m+k-1}})_{i=1}^{\infty}))$
  \\
  $\leqslant~\frac{\beta\left(d(x_{n+k-1},x_{m+k-1})\right)}{2}(d(T((x_{n+i-1,\widehat{n+k-1}})_{i=1}^{\infty}),x_{n+k-1})+d(T((x_{m+i-1,\widehat{m+k-1}})_{i=1}^{\infty}),x_{m+k-1}))$
  \\
  $<\frac{1}{2}(d(x_{n+k-1},x_{n+k})+d(x_{m+k-1},x_{m+k}))$.

  \par For sufficiently large enough $n,m\in\mathbb{N}$
  \[d(x_{n+k},x_{m+k})<\varepsilon\]
  for a fixed $\varepsilon>0$.

  Then the sequence $\{x_{n+k}\}_{n\in\mathbb{N}}$ is a Cauchy sequence.

  Since $(X,d)$ is a complete metric space, there exists a $u\in~X$ such that,
  \begin{equation}\label{2}
   \lim\limits_{n\to\infty}x_{n+k}=u.
  \end{equation}

  \par We claim that $T((u)_{i=1}^{\infty})=u$.

\par If we suppose the contrary we have $d(u,T((u)_{i=1}^{\infty}))>0$.

  Then, by (\ref{1}) and (\ref{2}), for an arbitrary $\varepsilon>0$ and a sufficiently large $n$ we get
  \begin{align*}
d\left(u,T\left((u)_{i=1}^{\infty}\right)\right)&\leqslant d\left(x_{n+k},u\right)+d(x_{n+k},T((u)_{i=1}^{\infty}))\\ &=d\left(x_{n+k},u\right)+d\left(T\left((x_{n+i-1,\widehat{n+k-1}})_{i=1}^{\infty}\right),T((u)_{i=1}^{\infty})\right)
			\\ & \leqslant d\left(x_{n+k},u\right)+\frac{\beta\left(d(u,x_{n+k-1})\right)}{2}(d(T((x_{n+i-1,\widehat{n+k-1}})_{i=1}^{\infty}),x_{n+k-1})+d(T((u)_{i=1}^{\infty}),u))
			\\
			& \leqslant \frac{1}{2}~d\left(u,T((u)_{i=1}^{\infty})\right)+\frac{\varepsilon}{2}+\frac{\varepsilon}{2}
 \leqslant \frac{1}{2}~d(u,T((u)_{i=1}^{\infty})+\varepsilon
\\
			\hspace{-25pt} &  \implies\frac{1}{2}~d\left(u,T((u)_{i=1}^{\infty})\right)\leqslant\varepsilon. ~\text{ This is a contradiction.}
  \end{align*}

  Hence, $d(T((u)_{i=1}^{\infty}),u)=0\implies~T((u)_{i=1}^{\infty})=u.$

  Then the conclusion follows.
\end{flushleft}
\end{proof}

Next, we will provide a new result for multivalued proper extension of Theorem \ref{2.5} which is also a generalization of Kannan (Theorem \ref{2.3}) as a result of Theorem \ref{3.2}.

\begin{corollary}\label{3.3}
Let $(X,d)$ be a complete metric space and $T$ be a Kannan-Geraghty $H_k$ contraction. Then T has a fixed point and every $k$-Picard sequence for $T$ converges to $u$.
\end{corollary}

\begin{proof}
Let us choose $x_0, \ldots, x_{k}\in~X$. We define a $k$-Picard sequence by
\begin{equation}
  x_{n+k}:=T(x_n,\ldots,x_{n+k-1}).
\end{equation}

If we follow same steps as in the proof of Theorem \ref{3.2} we get the required fixed point.
\end{proof}

\begin{remark}\label{3.4} For $k=1$ we get Theorem \ref{2.5} in \cite{Kg} which proves that Corollary \ref{3.3} is a proper generalization of Theorem \ref{2.5}.
\end{remark}
\begin{theorem}\label{3.5}
Let $\left(X,d\right)$ be a complete metric space and $T:\prod_{l=1}^{\infty} X \rightarrow X$ be an extended Fisher-Geraghty $H_k$ contraction for some $ k\in\mathbb{N}$. Then there exists $u\in X$ such that $T((u)_{i=1}^{\infty})=u$ and for any $x_1, x_2, \ldots x_k$ the infinite $k$-Picard sequence converges to $u$.
\end{theorem}

\begin{proof}
  Let $x_1, \ldots, x_k\in X$.  For all $n\in\mathbb{N} $ we define infinite $k$- Picard sequence by
		\begin{equation}\label{s}
		x_{n+k}:=T\left((x_{n+i-1,\widehat{n+k-1}})_{i=1}^{\infty}\right).
		\end{equation}

		We claim $\lim\limits_{n\to\infty}d(x_{n+k},x_{n+k+1})=0$.

Then we have
		\begin{flushleft}
			$d\left(x_{n+k+1},x_{n+k+2}\right)=d\left(T\left((x_{n+i,\widehat{n+k}})_{i=1}^{\infty}\right),T\left((x_{n+i+1,\widehat{n+1+k}})_{i=1}^{\infty}\right)\right)$
			\\
			$\leqslant\frac{\beta\left(d(x_{n+k},x_{n+k+1})\right)}{2}(d(T\left((x_{n+i,\widehat{n+k}})_{i=1}^{\infty}\right),x_{n+k+1})+d(T\left((x_{n+i,\widehat{n+k+1}})_{i=1}^{\infty}\right),x_{n+k}))$
			\\
			$<~\frac{1}{2}(d(x_{n+k+1},x_{n+k+1})+d(x_{n+k},x_{n+k+2}))$
			\\
			$<~\frac{1}{2}(d(x_{n+k},x_{n+k+1})+d(x_{n+k+1},x_{n+k+2}))\implies d(x_{n+k+1},x_{n+k+2})<d(x_{n+k},x_{n+k+1})$
		\end{flushleft}

Hence our sequence is a monotone decreasing sequence bounded below. This implies thrat there exsists $r\geq~0$ such that
\begin{equation}\label{3}
  \lim\limits_{n\to~0}d(x_{n+k},x_{n+k+1})=r.
\end{equation}

We claim that $r=0$.

If we suppose the contrary and $r>0$ we get
\begin{flushleft}
  $\lim\limits_{n\to\infty}d(x_{n+k+1},x_{n+k+2})$

  $\leqslant\frac{\lim\limits_{n\to\infty}\beta\left(d(x_{n+k},x_{n+k+1})\right)}{2}(\lim\limits_{n\to\infty}d(x_{n+k},x_{n+k+1})+\lim\limits_{n\to\infty}d(x_{n+k+1},x_{n+k+2}))$

  $\implies\frac{2\lim\limits_{n\to\infty}d(x_{n+k+1},x_{n+k+2})}{\lim\limits_{n\to\infty}d(x_{n+k},x_{n+k+1})+\lim\limits_{n\to\infty}d(x_{n+k+1},x_{n+k+2})}\leqslant~\lim\limits_{n\to\infty}\beta\left(d(x_{n+k},x_{n+k+1})\right)$

  Using (\ref{3}) we obtain

  $\frac{2r}{2r}\leqslant~\lim\limits_{n\to\infty}\beta\left(d(x_{n+k},x_{n+k+1})\right)\implies\lim\limits_{n\to\infty}\beta\left(d(x_{n+k},x_{n+k+1})\right)\geqslant~1 $
\end{flushleft}

\par Since $\beta\in\mathbb{G}$ we have	$\lim\limits_{n\to\infty}\beta\left(d(x_{n+k},x_{n+k+1})\right)\leqslant 1$.
\par From the well know in the related literature, Sandwich Theorem we get
\begin{equation}\label{l}
 \lim\limits_{n\to\infty}\beta\left(d(x_{n+k},x_{n+k+1})\right)=1 \implies\lim\limits_{n\to\infty}d(x_{n+k},x_{n+k+1})=0.
\end{equation}

\begin{flushleft}

 We claim that
			$\left\lbrace x_{n+k}\right\rbrace_{n\in\mathbb{N}}$ is a Cauchy and we want to prove this by contradiction. Then, using the contrary there exists $\varepsilon>0$ such that we can find subsequences $\left\lbrace x_{m\left(p\right)+k}\right\rbrace_{p\in\mathbb{N}}$,$\left\lbrace x_{n\left(p\right)+k}\right\rbrace_{p\in\mathbb{N}}$ with $m\left(p\right)>n\left(p\right)>p$ such that for every $p$ we have
			\begin{equation}\label{e2}
			d\left(x_{m\left(p\right)+k},x_{n\left(p\right)+k}\right)\geqslant\varepsilon.
			\end{equation}
		
\par Moreover, correseponding to each $n(p)$ we can choose least of such $m\left(p\right)$ satisfying (\ref{e2}) so that,
			\begin{equation}\label{e3}
			d\left(x_{m\left(p\right)+k-1},x_{n\left(p\right)+k}\right)<\varepsilon.
			\end{equation}

			Using (\ref{e2}), (\ref{e3}) and the triangle inequality we get
			\begin{equation}\label{e4}
			d\left(x_{m\left(p\right)+k-1},x_{n\left(p\right)+k-1}\right)\leqslant d\left(x_{m\left(p\right)+k-1},x_{n\left(p\right)+k}\right)+d\left(x_{n\left(p\right)+k-1},x_{n\left(p\right)+k}\right)
			\end{equation}
			\begin{center}
				$<\varepsilon+d\left(x_{n\left(p\right)+k-1},x_{n\left(p\right)+k}\right)$
			\end{center}
			and
			\begin{align*}\label{e5}
			\varepsilon &\leqslant d\left(x_{n\left(p\right)+k},x_{m\left(p\right)+k}\right)\\
			&\leqslant d\left(x_{n\left(p\right)+k},x_{n\left(p\right)+k-1}\right)+d\left(x_{n\left(p\right)+k-1},x_{m\left(p\right)+k-1}\right)+d\left(x_{m\left(p\right)+k-1},x_{m\left(p\right)+k}\right).
			\end{align*}

			Letting $p\to\infty$ in (\ref{e4}) and (\ref{l}) we get
			\begin{equation}\label{e6}
			\lim\limits_{p\to\infty}d\left(x_{m\left(p\right)+k-1},x_{n\left(p\right)+k-1}\right)=\varepsilon.
			\end{equation}

\begin{flushleft}
    $d(x_{n(p)+k},x_{m(p)+k})=d(T((x_{n(p)+i-1,\widehat{n(p)+k-1}})_{i=1}^{\infty}),T((x_{m(p)+i-1,\widehat{m(p)+k-1}})_{i=1}^{\infty}))$
    \\
    $\leqslant\frac{\beta(d(x_{n(p)+k-1},x_{m(p)+k-1}))}{2}(d(T((x_{n(p)+i-1,\widehat{n(p)+k-1}})_{i=1}^{\infty}),x_{m(p)+k-1})$
    \\
    $+d(T((x_{m(p)+i-1,\widehat{m(p)+k-1}})_{i=1}^{\infty}),x_{n(p)+k-1}))$\
    \\$<\frac{1}{2}(d(x_{n(p)+k},x_{m(p)+k-1})+d(x_{m(p)+k},x_{n(p)+k-1}))$
    \\$\leqslant\frac{1}{2}(\varepsilon+d(x_{m(p)+k},x_{m(p)+k-1})+d(x_{m(p)+k-1},x_{n(p)+k})+d(x_{n(p)+k},x_{n(p)+k-1}))$
    \\$\leqslant\frac{1}{2}(2\varepsilon+d(x_{m(p)+k},x_{m(p)+k-1})+d(x_{n(p)+k},x_{n(p)+k-1}))$
\end{flushleft}

    Then from (\ref{e6}) and (\ref{l}) we get $\varepsilon<\varepsilon$. So this is a contradiction.

    Then
    \begin{equation}
     \{x_{n+k}\}_{n\in\mathbb{N}}~\text{ is a Cauchy sequence.}
    \end{equation}

  Since $(X,d)$ is complete there exists a $u\in~X$ such that
  \begin{equation}\label{4}
   \lim\limits_{n\to\infty}x_{n+k}=u.
  \end{equation}
 \par  We claim that $T((u)_{i=1}^{\infty})=u$.

 \par If we suppose the contrary we have $d(u,T((u)_{i=1}^{\infty}))>0$. Then, by (\ref{3}) and (\ref{4}) for arbitrary an $\varepsilon>0$ and a sufficiently large $n$ we get
  \begin{align*}
&d\left(u,T\left((u)_{i=1}^{\infty}\right)\right)\leqslant d\left(x_{n+k},u\right)+d(x_{n+k},T((u)_{i=1}^{\infty}))\\ &=d\left(x_{n+k},u\right)+d\left(T\left((x_{n+i-1,\widehat{n+k-1}})_{i=1}^{\infty}\right),T((u)_{i=1}^{\infty})\right)
			\\ & \leqslant d\left(x_{n+k},u\right)+\frac{\beta\left(d(u,x_{n+k-1})\right)}{2}(d(T((x_{n+i-1,\widehat{n+k-1}})_{i=1})^{\infty}),u)+d(T((u)_{i=1}^{\infty}),x_{n+k-1}))
			\\
			& \leqslant d\left(x_{n+k},u\right)+\frac{\beta\left(d(u,x_{n+k-1})\right)}{2}(d(T((x_{n+i-1,\widehat{n+k-1}})_{i=1})^{\infty}),u)+d(T((u)_{i=1}^{\infty}),u)+d(u,x_{n+k-1}))
			\\
			& \leqslant \frac{1}{2}~d\left(u,T((u)_{i=1}^{\infty})\right)+\frac{\varepsilon}{3}+\frac{\varepsilon}{6}
			\\
			& \leqslant \frac{1}{2}~d(u,T((u)_{i=1}^{\infty})+\frac{1}{2}\varepsilon\\
			\hspace{-25pt} &  \implies d\left(u,T((u)_{i=1}^{\infty})\right)\leqslant\varepsilon.
  \end{align*}

 \par Hence $d(T((u)_{i=1}^{\infty}),u)=0$ which implies $T((u)_{i=1}^{\infty})=u.$ Then the conclusion follows.

\end{flushleft}
\end{proof}

\begin{corollary}\label{3.6}
Let $(X,d)$ be a complete metric space and $T$ be a Fisher-Geraghty $H_k$ contraction. Then $T$  has a fixed point.
\end{corollary}
\begin{proof}

Choosen any $x_0, \ldots, x_{k}\in~X$ we define
\begin{equation}
  x_{n+k}:=T(x_n,\ldots,x_{n+k-1}).
\end{equation}
Using same steps as in the proof of Theorem \ref{3.5} we get the conclusion.
\end{proof}

\begin{remark}\label{3.7} For $k=1$ we get a new type of extension of Theorem \ref{2.3} proved in \cite{F} and for any $k$, Corollary \ref{3.6} also gives the multidimensional extension of the same theorem stated in \cite{F}.
\end{remark}
\begin{remark} The Banach Fixed Point Theorem(\cite{Ban}), Theorem \ref{2.1} , Theorem \ref{2.2} , Theorem \ref{2.3}, Theorem \ref{2.5} and Theorem \ref{2.6}  are all applicable only in complete metric space as well as the theorems described in \cite{Kho} but  Theorem \ref{3.1}, Theorem \ref{3.2}, Theorem \ref{3.5} which we have proved are talking about the space $\prod_{i=1}^{\infty}X$( where $X$ is a complete metric space) which is indeed metrizable but might not be complete.
\end{remark}

\section{Applications to integral equations}

\quad In \cite{n-dim} H Singh et al. introduced the notion of multi-dimensional Fredholm integral equation taking into account $n\in\mathbb{N}$ as dimension of the equation.
$$u(x_1,...,x_n)=f(x_1,...,x_n)+\underbrace{\int\limits_{0}^{1}\int\limits_{0}^{1}...\int\limits_{0}^{1}}_{n-times}K(x_1,...,x_n,y_1,...,y_n)u(y_1,...,y_n)dy_1...dy_n,$$
where $(x_1,...,x_n)\in D=\underbrace{([0,1]\times[0,1]\times...\times [0,1])}_{n-times}$ and $f(x_1,...,x_n), K(x_1,...,x_n,y_1,...,y_n)$ are known continuous functions defined on $D$ and $D^2$ respectively and $u(x_1,...,x_n)$ is a unknown function.

The Fredholm integral equations play an important role in modelling of physics phenomena described by two or three dimensions. In the same way, they have applications in astrophysics models thinking of the four dimensions of a neutron star or a black hole.

Thinking of this aspects, if we extend to infinity the dimension "n" of the previous multi-dimensional integral equation we introduce a new notion, {\it the infinite dimensional Fredholm integral equation}, as follows.

\begin{equation}\label{w1}
u\left((t_{i,\widehat{k}})_{i=1}^{\infty}\right)=f\left((t_{i,\widehat{k}})_{i=1}^{\infty}\right)+\int\limits_{0}^{1}...\int\limits_{0}^{1}...K\left( (t_{i,\widehat{k}})_{i=1}^{\infty},(s_{i,\widehat{k}})_{i=1}^{\infty}\right)u\left((s_{i,\widehat{k}})_{i=1}^{\infty}\right)ds_1...ds_n... .
\end{equation}

where $f\left((t_{i,\widehat{k}})_{i=1}^{\infty}\right):[0,1]\to\mathbb{R}$ and $K:[0,1]\times\mathbb{R}\to\mathbb{R}$ are two known continuous functions and $u\left((t_{i,\widehat{k}})_{i=1}^{\infty}\right)$ is an unknown function.

Further, let us give our first application of the main results of this paper proving the existence of a solution of infinite dimensional Fredholm integral equation (\ref{w1}).

\begin{theorem}

Let $X=C([0,1],\mathbb{R})$  the set of real continuous functions on $[0,1]$ and let $d:X\times X\to\mathbb{R}_{+}$ given by
\begin{equation}\label{w2}
d\left((u_{i,\widehat{k}})_{i=1}^{\infty},T(u_{i,\widehat{k}})_{i=1}^{\infty}\right)=\sup\limits_{t\in[0,1]}\mid (u_{i,\widehat{k}})_{i=1}^{\infty}-T(u_{i,\widehat{k}})_{i=1}^{\infty}\mid.\end{equation}
Define $T:\prod\limits_{i=1}^{\infty}X\to X$ by
\begin{equation}\label{w3}
Tu\left((t_{i,\widehat{k}})_{i=1}^{\infty}\right)=f\left((t_{i,\widehat{k}})_{i=1}^{\infty}\right)+\int\limits_{0}^{1}...\int\limits_{0}^{1}...K\left( (t_{i,\widehat{k}})_{i=1}^{\infty},(s_{i,\widehat{k}})_{i=1}^{\infty}\right)u\left((s_{i,\widehat{k}})_{i=1}^{\infty}\right)ds_1...ds_n... .
\end{equation}

Assume the following holds:

\begin{enumerate}
  \item[$(i)$]$\left| u\left((t_{i,\widehat{k}})_{i=1}^{\infty}\right)-Tu\left((t_{i,\widehat{k}})_{i=1}^{\infty}\right)\right|\leq\frac{1}{2}\left| d\left(T\left((u_{i,\widehat{k}})_{i=1}^{\infty}\right), u_k\right)+d\left(T\left((u_{i+1,\widehat{k+1}})_{i=1}^{\infty}\right),u_{k+1}\right)\right|$;
  \item[$(ii)$] there exists a constant $\delta\in(0,1)$ such that $K\left((t_{i,\widehat{k}})_{i=1}^{\infty},(s_{i,\widehat{k}})_{i=1}^{\infty}\right)<\delta$;
  \item[$(iii)$]let $\beta:[0,\infty)\to[0,1)$ be defined as $\beta\left(z\right)=\frac{1}{\gamma}<1$ for every $\gamma>0$ and $z\in \mathbb{R}_{+}$.
\end{enumerate}

Then the infinite dimensional Fredholm integral equation (\ref{w1}) has a solution.

\end{theorem}

\begin{proof}

We can prove the existence of a solution of infinite dimensional Fredholm integral equation if we show that the operator $T$ defined by (\ref{w3}) has a fixed point.

We can easy remark that the space $X=(C[0,1],\mathbb{R})$ endowed with the metric $d$ defined by relation (\ref{w2}) form a complete metric space.
Then we shall show that all the hypothesis of Theorem \ref{3.2} are verified.

We have the following estimation
\begin{align*}
  &\left| Tu\left((t_{i,\widehat{k}})_{i=1}^{\infty}\right)-T^2u\left((t_{i,\widehat{k}})_{i=1}^{\infty}\right)\right| =|\int\limits_{0}^{1}...\int\limits_{0}^{1}...K\left( (t_{i,\widehat{k}})_{i=1}^{\infty},(s_{i,\widehat{k}})_{i=1}^{\infty}\right)u\left((s_{i,\widehat{k}})_{i=1}^{\infty}\right)ds_1...ds_n...\\
  & - \int\limits_{0}^{1}...\int\limits_{0}^{1}...K\left( (t_{i,\widehat{k}})_{i=1}^{\infty},(s_{i,\widehat{k}})_{i=1}^{\infty}\right)Tu\left((s_{i,\widehat{k}})_{i=1}^{\infty}\right)ds_1...ds_n...| \\
  & \leq \int\limits_{0}^{1}...\int\limits_{0}^{1}...K\left( (t_{i,\widehat{k}})_{i=1}^{\infty},(s_{i,\widehat{k}})_{i=1}^{\infty}\right)\left|u\left((s_{i,\widehat{k}})_{i=1}^{\infty}\right)-Tu\left((s_{i,\widehat{k}})_{i=1}^{\infty}\right)\right|ds_1...ds_n...\\
  & \leq \frac{\delta}{2\gamma}\left| d\left(T\left((u_{i,\widehat{k}})_{i=1}^{\infty}\right), u_k\right)+d\left(T\left((u_{i+1,\widehat{k+1}})_{i=1}^{\infty}\right),u_{k+1}\right)\right|\int\limits_{0}^{1}...\int\limits_{0}^{1}...ds_1...ds_n...\\
  &=\frac{\delta}{2\gamma}\left| d\left(T\left((u_{i,\widehat{k}})_{i=1}^{\infty}\right), u_k\right)+d\left(T\left((u_{i+1,\widehat{k+1}})_{i=1}^{\infty}\right),u_{k+1}\right)\right|.
\end{align*}

Taking supremum on both sides we get

$$\sup\limits_{t\in[0,1]}\left| Tu\left((t_{i,\widehat{k}})_{i=1}^{\infty}\right)-T^2u\left((t_{i,\widehat{k}})_{i=1}^{\infty}\right)\right|$$
$$\leq \sup\limits_{t\in[0,1]}\frac{\delta }{2\gamma}\left| d\left(T\left((u_{i,\widehat{k}})_{i=1}^{\infty}\right), u_k\right)+d\left(T\left((u_{i+1,\widehat{k+1}})_{i=1}^{\infty}\right),u_{k+1}\right)\right|.$$

Then, for $z=d(u_k,u_{k+1})$ and $\delta\in(0,1)$ we obtain

\begin{align*}
  d\left(Tu\left((t_{i,\widehat{k}})_{i=1}^{\infty}\right),T^2u\left((t_{i,\widehat{k}})_{i=1}^{\infty}\right)\right)\\&\hspace{-70pt}\leq \frac{\delta\beta(d(u_k,u_{k+1}))}{2}\left(d\left(T\left((u_{i,\widehat{k}})_{i=1}^{\infty}\right), u_k\right)+d\left(T\left((u_{i+1,\widehat{k+1}})_{i=1}^{\infty}\right),u_{k+1}\right) \right)\\
  &\hspace{-70pt}\leq \frac{\beta(d(u_k,u_{k+1}))}{2}\left(d\left(T\left((u_{i,\widehat{k}})_{i=1}^{\infty}\right), u_k\right)+d\left(T\left((u_{i+1,\widehat{k+1}})_{i=1}^{\infty}\right),u_{k+1}\right) \right)
\end{align*}

In conclusion, all the hypothesis of  Theorem \ref{3.2} are accomplished. Then the operator $T$ has a fixed point, which means the infinite dimensional Fredholm integral equation (\ref{w1}) has a solution.
\end{proof}

The following application involve another type of integral equations, Urysohn type integral equations. We extend the known cases of this type of integral equations to {\it infinite dimensional Urysohn integral equation}.
\begin{equation}\label{urh}
  u\left((t_{i,\widehat{k}})_{i=1}^{\infty}\right)=f\left((t_{i,\widehat{k}})_{i=1}^{\infty}\right)+\int\limits_{a}^{b}...\int\limits_{a}^{b}... P\left( (t_{i,\widehat{k}})_{i=1}^{\infty},(s_{i,\widehat{k}})_{i=1}^{\infty},u\left((s_{i,\widehat{k}})_{i=1}^{\infty}\right)\right)ds_1...ds_n... .
\end{equation}
where $f\left((t_{i,\widehat{k}})_{i=1}^{\infty}\right):[a,b]\to\mathbb{R}$ and $K:[a,b]\times[a,b]\times\mathbb{R}\to\mathbb{R}$ are two known continuous functions and $u\left((t_{i,\widehat{k}})_{i=1}^{\infty}\right)$ is an unknown function.

For this new type of Urysohn integral equation let us give the following result.

\begin{theorem}
Let $X=C([a,b],\mathbb{R})$  the set of real continuous functions on $[a,b]$ and let $d:X\times X\to\mathbb{R}_{+}$ given by
\begin{equation}\label{urh1}
d\left((u_{i,\widehat{k}})_{i=1}^{\infty},T(u_{i,\widehat{k}})_{i=1}^{\infty}\right)=\max\limits_{t\in[a,b]}\mid (u_{i,\widehat{k}})_{i=1}^{\infty}-T(u_{i,\widehat{k}})_{i=1}^{\infty}\mid\cdot e^{-sin^{-1}\alpha},\end{equation}
with $\alpha>0$.

Define $T:\prod\limits_{i=1}^{\infty}X\to X$ by
\begin{equation}\label{urh2}
Tu\left((t_{i,\widehat{k}})_{i=1}^{\infty}\right)=f\left((t_{i,\widehat{k}})_{i=1}^{\infty}\right)+\int\limits_{a}^{b}...\int\limits_{a}^{b}...P\left( (t_{i,\widehat{k}})_{i=1}^{\infty},(s_{i,\widehat{k}})_{i=1}^{\infty},u\left((s_{i,\widehat{k}})_{i=1}^{\infty}\right)\right)ds_1...ds_n... .
\end{equation}

Assume the following holds:

\begin{enumerate}
  \item[$(i)$]$\left| P\left((t_{i,\widehat{k}})_{i=1}^{\infty},(s_{i,\widehat{k}})_{i=1}^{\infty},u\left((s_{i,\widehat{k}})_{i=1}^{\infty}\right)\right)-P\left((t_{i,\widehat{k}})_{i=1}^{\infty},(s_{i,\widehat{k}})_{i=1}^{\infty},Tu\left((s_{i,\widehat{k}})_{i=1}^{\infty}\right)\right)\right|\\ \leq\frac{1}{\tau}\left| (u_{i,\widehat{k}})_{i=1}^{\infty}-T(u_{i,\widehat{k}})_{i=1}^{\infty}\right|$;
   \item[$(ii)$]let $\beta:[0,\infty)\to[0,1)$ be defined as $\beta\left(z\right)=\frac{1}{\tau}<1$ for every $\tau>0$ and $z\in X$.
\end{enumerate}

Then the infinite dimensional Urysohn integral equation (\ref{urh}) has a solution.

\end{theorem}

\begin{proof}

It is easy to check the space $X=(C[0,1],\mathbb{R})$ endowed with the metric $d$ defined by relation (\ref{urh1}) is a complete metric space.

To prove the existence of a solution of infinite dimensional Urysohn integral equation we shall show that the operator $T$ defined by (\ref{urh2}) has a fixed point.

We have the following estimation
\begin{align*}
  &\left| Tu\left((t_{i,\widehat{k}})_{i=1}^{\infty}\right)-T^2u\left((t_{i,\widehat{k}})_{i=1}^{\infty}\right)\right| =|\int\limits_{a}^{b}...\int\limits_{a}^{b}...P\left( (t_{i,\widehat{k}})_{i=1}^{\infty},(s_{i,\widehat{k}})_{i=1}^{\infty},u\left((s_{i,\widehat{k}})_{i=1}^{\infty}\right)\right)ds_1...ds_n...\\
  & - \int\limits_{a}^{b}...\int\limits_{a}^{b}...P\left( (t_{i,\widehat{k}})_{i=1}^{\infty},(s_{i,\widehat{k}})_{i=1}^{\infty},Tu\left((s_{i,\widehat{k}})_{i=1}^{\infty}\right)\right)ds_1...ds_n...| \\
  &\leq \int\limits_{a}^{b}...\int\limits_{a}^{b}...\left| P\left((t_{i,\widehat{k}})_{i=1}^{\infty},(s_{i,\widehat{k}})_{i=1}^{\infty},u\left((s_{i,\widehat{k}})_{i=1}^{\infty}\right)\right)-P\left((t_{i,\widehat{k}})_{i=1}^{\infty},(s_{i,\widehat{k}})_{i=1}^{\infty},Tu\left((s_{i,\widehat{k}})_{i=1}^{\infty}\right)\right)\right|ds_1...ds_n...\\
  & \leq \frac{1}{\tau}\int\limits_{a}^{b}...\int\limits_{a}^{b}... \left| (u_{i,\widehat{k}})_{i=1}^{\infty}-T(u_{i,\widehat{k}})_{i=1}^{\infty}\right|ds_1...ds_n...\\
  & \leq \frac{e^{sin^{-1}\alpha}}{\tau}\left|(u_{i,\widehat{k}})_{i=1}^{\infty}-T(u_{i,\widehat{k}})_{i=1}^{\infty}\right|e^{-sin^{-1}\alpha}\int\limits_{a}^{b}...\int\limits_{a}^{b}...ds_1...ds_n...\\
  &\leq \frac{e^{sin^{-1}\alpha}}{\tau}\left|(u_{i,\widehat{k}})_{i=1}^{\infty}-T(u_{i,\widehat{k}})_{i=1}^{\infty}\right|e^{-sin^{-1}\alpha}.
   \end{align*}

Applying maximum on both sides we get
$$\max\limits_{t\in[a,b]}\left| Tu\left((t_{i,\widehat{k}})_{i=1}^{\infty}\right)-T^2u\left((t_{i,\widehat{k}})_{i=1}^{\infty}\right)\right|e^{-sin^{-1}\alpha}\leq \frac{1}{\tau}\max\limits_{t\in[a,b]} \left|(u_{i,\widehat{k}})_{i=1}^{\infty}-T(u_{i,\widehat{k}})_{i=1}^{\infty}\right|e^{-sin^{-1}\alpha}.$$

For $z=M_{k}\left((u_{i,\widehat{k}})_{i=1}^{\infty},(u_{i+1,\widehat{k+1}})_{i=1}^{\infty}\right)$ we get
\begin{align*}
d\left(Tu\left((t_{i,\widehat{k}})_{i=1}^{\infty}\right),T^2u\left((t_{i,\widehat{k}})_{i=1}^{\infty}\right)\right)&\leq \beta\left(M_{k}\left((u_{i,\widehat{k}})_{i=1}^{\infty},(u_{i+1,\widehat{k+1}})_{i=1}^{\infty}\right)\right)d\left((u_{i,\widehat{k}})_{i=1}^{\infty},T(u_{i,\widehat{k}})_{i=1}^{\infty}\right)\\
&\leq \beta\left(M_{k}\left((u_{i,\widehat{k}})_{i=1}^{\infty},(u_{i+1,\widehat{k+1}})_{i=1}^{\infty}\right)\right)M_{k}\left((u_{i,\widehat{k}})_{i=1}^{\infty},(u_{i+1,\widehat{k+1}})_{i=1}^{\infty}\right)
\end{align*}

Then all the conditions of Theorem \ref{3.1} are satisfied; result that the operator $T$ has a fixed point.

Then the infinite dimensional Urysohn type integral equation (\ref{urh}) has a solution.
\end{proof}

\section{Conclusion}

\quad Here in this article we have shown some metric combinatorial arguments that can be profitably applied to extend Geraghty's theorem, Kannan-Geraghty theorem and, Fisher's theorem as well as a particular extension for the infinite cases. This idea could also motivate to find a  similar infinite extension of the other contractions and its related operators to find fixed points or coincidence points. We also defined two new types of infinite integral equations which would encourage the study of the infinite integral equations and their solutions. From the discussion we have some interesting problems as follows.

\begin{enumerate}
 \item Can a similar concept be extended to the notion of coincidence points?

 \item We have discussed the fixed point theorem for the contraction from countable products of the complete metric space $X$ (with respect to uniform metric topology) to the space $X$. Can it further be extended to  an uncountable product of complete metric space which is not necessarily a metric space?
\end{enumerate}

%\noindent \textbf{Conflicts of Interest}
%Authors declare that there is no conflicts of interest regarding the publication of this paper.\\

\noindent \textbf{Acknowledgements} The  authors would like to thank Suprokash Hazra who has helped by his valuable comments in developing the key ideas of this paper.

\end{document}